\documentclass[11pt]{article}
\usepackage{graphicx}
\usepackage{amssymb}
\usepackage{epstopdf}
\usepackage{amsmath}
\usepackage{amsfonts}
\usepackage{amsthm}
\usepackage{latexsym}
\usepackage{url}
\usepackage{verbatim}
\usepackage{float}
\usepackage[all]{xy}
\def\Z{{\bf{Z}}}
\def\Q{{\bf{Q}}}
\def\C{{\bf{C}}}

\def\Hom{\mathrm{Hom}}

\title{Robert F. Coleman 1954 -- 2014}
\author{Matthew Baker, Barry Mazur, Ken Ribet}
\begin{document}

\maketitle

%\tableofcontents
%\part{Introduction}

\begin{quote} 
IN MEMORIAM:   
{\it
Robert Coleman
Professor of Mathematics,
UC Berkeley
%1954 -- 2014
}\end{quote}

\begin{figure}[H]
\centering
\includegraphics[width=.4\textwidth]{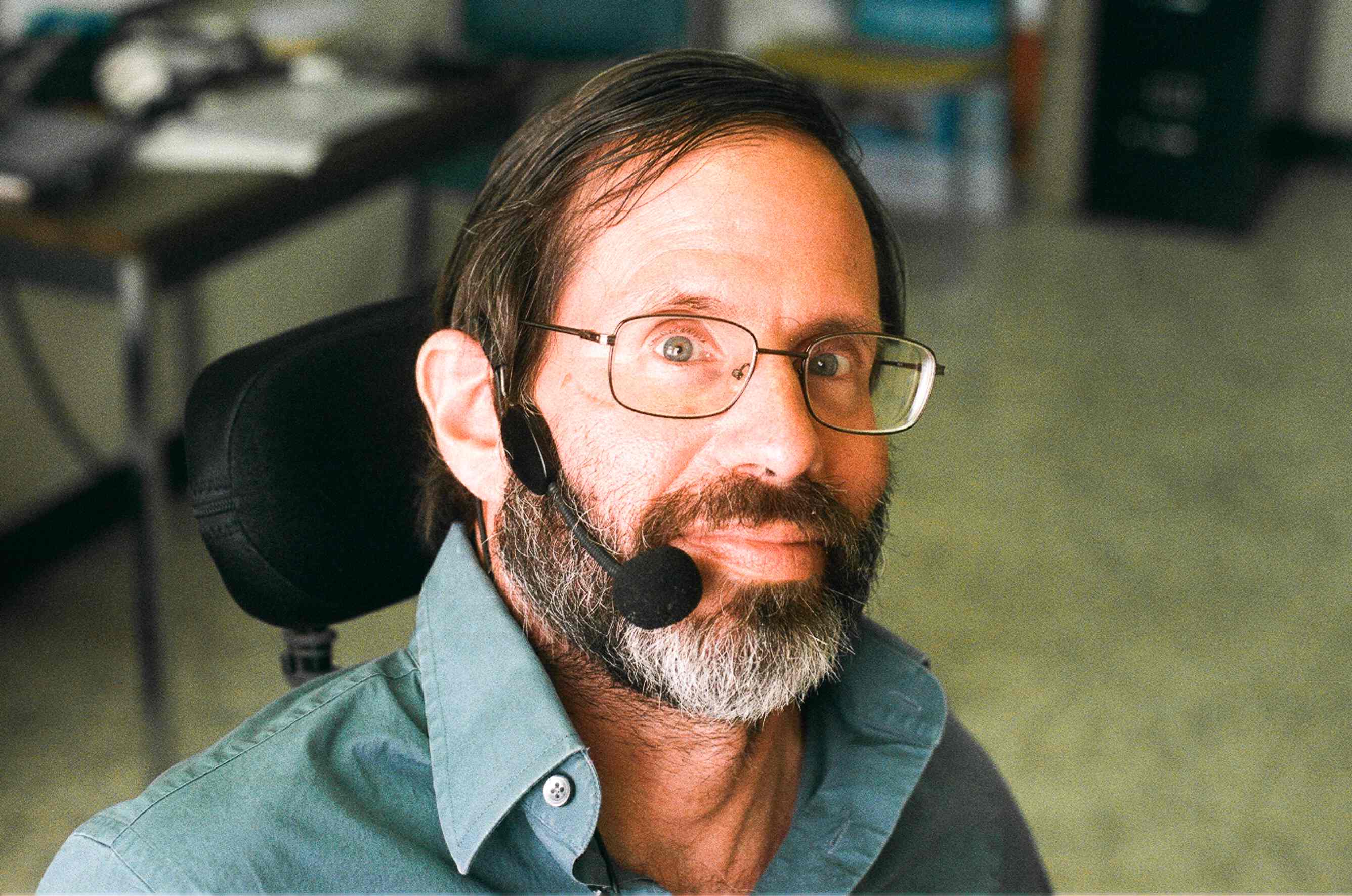}
\end{figure}

%  \hskip100pt \includegraphics[width=0.4\textwidth]{bishop.jpg}

\section{Biography}

Robert F. Coleman, an extraordinarily original and creative mathematician who has had a profound influence on modern number theory and arithmetic geometry, died of a sudden heart attack in El Cerrito, CA on the morning of March  24, 2014. He is survived by his wife Tessa, his sister Rosalind and brother Mark, his nephew Jeffrey and niece Elise, and his service dog Julep. The depth and importance of his mathematical ideas, his congeniality, the joy radiating from his playful disposition, and his sheer inexhaustible energy---all  this in the face of Multiple Sclerosis, a condition that did not deter him from full engagement with life---made Robert an inspiration to his friends, family, students, and colleagues.  Robert also worked toward making civic structures and laws more appropriate for people with disabilities, and his activism is yet another reason that Robert was so widely admired.  
% He will be dearly missed.

\medskip 

Robert was born on Nov. 22, 1954 in Glen Cove, NY.  He displayed an early talent for mathematics, winning an Intel Science Talent Search Award in 1972 as a high school student.  He earned a mathematics degree from Harvard University and subsequently completed Part III of the mathematical tripos at Cambridge, where he did research under the supervision of John Coates.  By the time he entered graduate school at Princeton, Robert had essentially already written his doctoral dissertation, but his formal thesis advisor was Kenkichi Iwasawa.  His dissertation, entitled ``Division Values in Local Fields,"  is considered a landmark contribution to local class field theory.  After completing his Ph.D., Robert returned to his alma mater, Harvard University, as a Benjamin Peirce Assistant Professor and Research Associate. He came to UC Berkeley as an Assistant Professor in 1983, and was promoted unusually quickly to Associate and then Full Professor.  He taught and did mathematics at Berkeley until his untimely death.  He had 12 Ph.D. students, published 63 research papers, and received numerous honors and recognition for his work, including a MacArthur ``Genius" Fellowship in 1987. 

\begin{figure}[H]
\centering
\includegraphics[width=.4\textwidth]{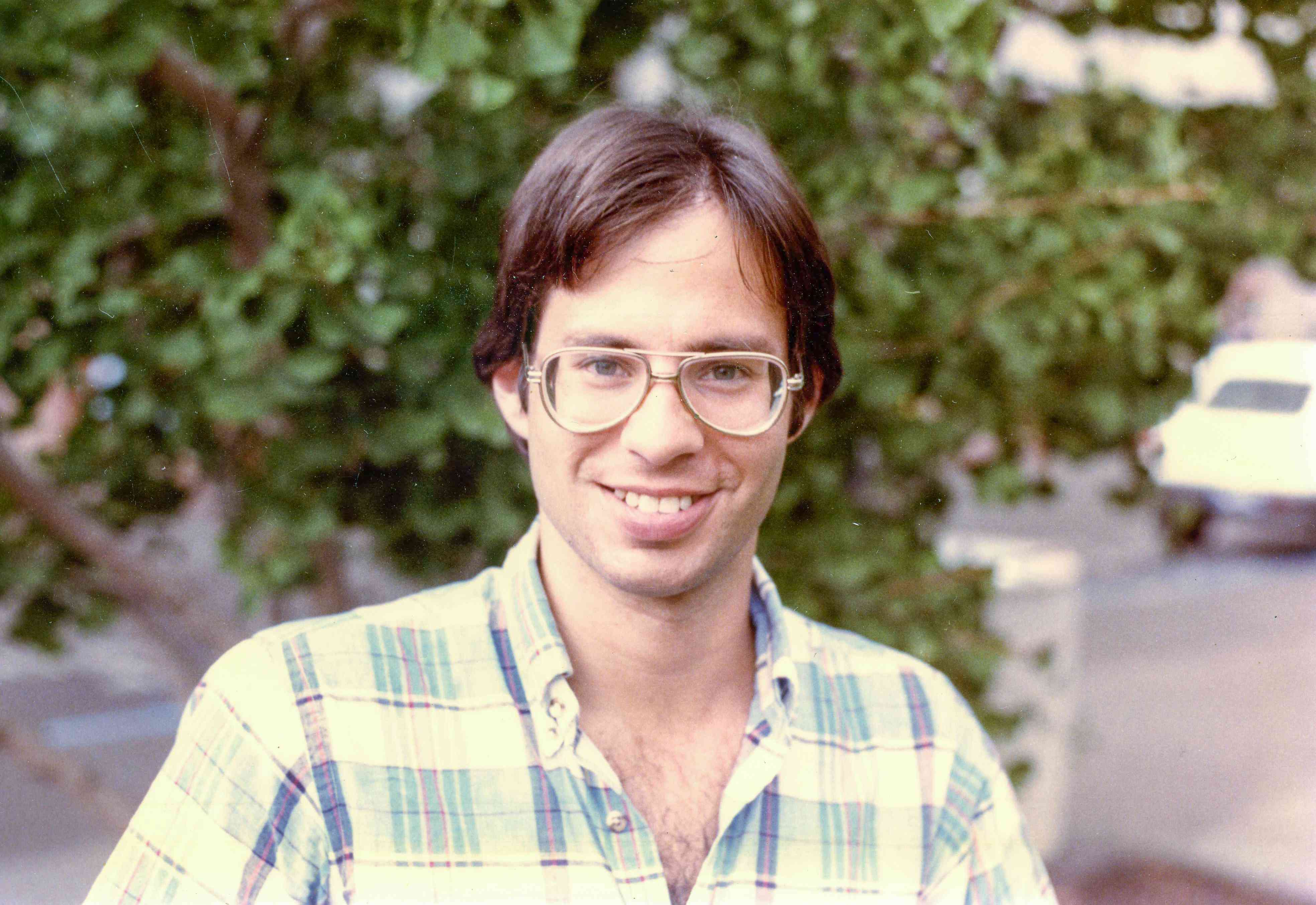}
\end{figure}

  % including 8 papers in the prestigious journal Inventiones Mathematicae and 5 in the Duke Mathematical Journal. 
  
   % He had an amazing intuition for everything p-adic.  Long before the invention of Berkovich spaces, Robert could somehow visualize paths and structures in p-adic geometry which no one else in the world saw as keenly or as profoundly
  
\section{Overview of Robert's mathematics}
  
Robert Coleman's primary mathematical love was number theory, with a particular interest in $p$-adic aspects of the subject\footnote{In the early twentieth century, Hensel and others introduced, for each prime number $p$, a field $\Q_p$ of $p$-adic numbers which serves as an ``arithmetic'' analogue of the real numbers.  There is also a larger algebraically closed field $\C_p$ corresponding to the complex numbers.  It was soon realized that not only do $p$-adic numbers have many important applications, but many facts from real and complex analysis and complex algebraic geometry have quite precise $p$-adic analogues.}, for which he had an amazing intuition.  
We will briefly summarize here some of Robert's key contributions to mathematics; these will be elaborated upon in subsequent sections.
% We'll give the flavor of this selection of his fundamental contributions to our subject, in the sections below.  
  \begin{itemize} 
  
     \item[]   {\bf Coleman maps:}  Classical analytic number theory, as developed in the works of Dirichlet and Dedekind, and in Minkowski's Geometry of Numbers, makes a deep connection between arithmetic and analysis in algebraic number theory, with the theory of {\it the regulator} as the vital bridge. 
In the $p$-adic analytic number theory of number fields, elliptic curves, and modular forms, ``Coleman maps'" provide the corresponding $p$-adic bridge. We will elaborate on this theme in Section~\ref{reg}.
% In close analogy, the theory of {\em $p$-adic regulators} is of fundamental importance for the $p$-adic-analytic number theory of number fields, elliptic curves, and modular forms. ``Coleman maps'" provide the corresponding $p$-adic bridge. {\bf Matt: The non-parallelism of these two sentences needs to be addressed, and we need a little more background history here.}   We will elaborate on this theme in section \ref{reg} below.
     
  \item[] {\bf Coleman integration:}  
% The $p$-adic fields were initially developed, by Hensel and others at the end of the nineteenth century and were immediately seen to be thorough-going analogues---in the role they played in arithmetic geometry--to the real and complex fields, ${\R}$ and ${\C}$. 
There is a serious obstacle to transferring one of our most useful tools in complex analysis---integration over paths---from ${\C}$ to $\C_p$.  Among other obstructions, one of the fundamental problems is the totally disconnected topology of the $p$-adic fields.  
Robert Coleman could somehow visualize paths and structures in $p$-adic geometry which no one else in the world saw as keenly or as profoundly{\footnote{This was long before the invention of Berkovich spaces, which provide a natural way to embed $\C_p$ in a path-connected space.}}. Thanks to Robert's vision, we can now integrate $p$-adic differential one-forms over paths in a manner analogous to the classical complex theory.  
Robert called the key new idea ``analytic continuation along Frobenius''.  In his 1985 Annals of Math paper ``Torsion Points on Curves and $p$-adic Abelian Integrals", he wrote: 
  \begin{quote} Rigid analysis was created to provide some coherence in an otherwise totally disconnected $p$-adic realm. Still, it is often left to Frobenius to quell the rebellious outer provinces.\end{quote}
    Applications of Coleman integration include the celebrated ``method of Chabauty--Coleman'' for finding rational points on curves, a new proof of the Manin-Mumford conjecture, and new results in $p$-adic Hodge theory.
    We will elaborate on these themes in Section~\ref{p-int}.
  
     \item[] {\bf Coleman Families:}  The theme of ``interpolation" of classical modular forms of integral weight to get ``continuous families'' of $p$-adic weights is one that might be traced back (with appropriate translation of language) to Ramanujan. The emergence of the interpolated ``$p$-adic modular eigenforms"  is one of the most important topics of contemporary number theory.  Following pioneering work of J.-P. Serre, Robert produced {\it $p$-adic families} of $p$-adic Banach spaces  of (overconvergent) $p$-adic modular forms, together with a completely continuous operator (the Atkin--Lehner operator) whose Fredholm spectrum gives one a good handle on the locus of eigenforms.  Robert's construction also opened the way for  the ``eigencurve'', about which we will elaborate in Section~\ref{eigen}.
   
 \end{itemize}

All three of these directions of Robert's research dovetail together, and remain thoroughly contemporary in their importance.  This point was illustrated in a comment made by Kevin Buzzard\footnote{Every personal reminiscence in this article whose source is not identified is drawn from the comments section in the first author's blog post \url{http://mattbakerblog.wordpress.com/2014/03/25/robert-f-coleman-1954-2014/}} on the day that he heard of Robert's death:
  
  \begin{quote} Two or three years ago, I remember noting one afternoon [after hearing three talks at Imperial College] that even though the topics were covering a lot of modern algebraic number theory ($(\phi,\Gamma)$-modules and their applications, $p$-adic modular forms, and computational methods for rational points on algebraic varieties), Robert's work had been a critical component of all three lectures I heard that day: the first was on the Coleman homomorphism, the second on Coleman families, and the third used Coleman's explicit Chabauty. A wonderful and genuinely coincidental moment of synchronicity which only underlined how much Robert has influenced modern number theory.\end{quote}

Robert made many other important contributions to number theory, too many to discuss in any depth in this article.  For example, he discovered and filled in a gap in Manin's proof of the Mordell conjecture over function fields and, with Voloch, verified some previously unchecked compatibilities in Dick Gross's theory of companion forms, thus completing Gross's influential work.  
% He found a totally new proof of the Manin--Mumford conjecture (originally proved by Raynaud), and used his approach (in joint work with Tamagawa and Tzermias) to explicitly compute the cuspidal torsion packet on Fermat curves.  
His work on semistable models for modular curves of arbitrary level was recently refined and applied to the Local Langlands Program by Jared Weinstein.

\section{The challenges of MS}

Robert Coleman was an avid tennis player and world traveler when he was struck with a severe case of Multiple Sclerosis in 1985.  According to his friend and collaborator Bill McCallum:

\begin{quote}[Robert] was on a visit to Japan and had a couple of strange incidents: being unable to hit the ball in tennis, stumbling on the stairs in the subway. The Japanese doctors didn't know what was going on. When he got back to Berkeley he was diagnosed with MS. From that point he had a shockingly rapid descent that caught his doctors by surprise; for many patients MS is a long slow decline. He went from perfectly healthy to the verge of death within a matter of weeks. His decline was arrested by some experimental and aggressive chemotherapy treatment\ldots
% . After that he recovered slightly, to the point where he could come home and get around with the aid of a wheel chair. 
But he never bounced back the way MS patients sometimes do.\end{quote}

Robert fought MS bravely and with a great sense of humor, and despite the severity of his illness he remained an active faculty member at UC Berkeley until his retirement in 2013.  Those that knew him were consistently amazed by his optimism and by the way he continued to travel despite the challenges of his MS.  His former Ph.D. student Harvey Stein called it ``truly amazing" how happily Robert lived his life, carrying on ``as if [his MS] wasn't even relevant".  
 His friend and colleague Jeremy Teitelbaum said:

\begin{quote} Of course Robert was brilliant, but what really impressed me about him was how relentless he was when it came to mathematics. He simply never gave up in the pursuit of an idea\dots. When his MS had made his hands shake and his voice slur, he brought that same determination to bear. I watched him spend hours typing slides for his lectures so that he could keep teaching---things that should have taken an hour would take half a day, but he did not give up.\end{quote}
 
Robert had a mischievous and impish sense of humor, and he surrounded himself with colorful and funny people.  For many years Robert's closest companion was his guide dog Bishop, who would join Robert everywhere.  Bishop eventually passed away and Robert found a new canine companion named Julep.  (Memorial fund donations to ``Paws With A Cause" can be made in honor of Robert's special relationship with his service dogs.) 
 
\medskip

\begin{figure}[H]
\centering
\includegraphics[width=.4\textwidth]{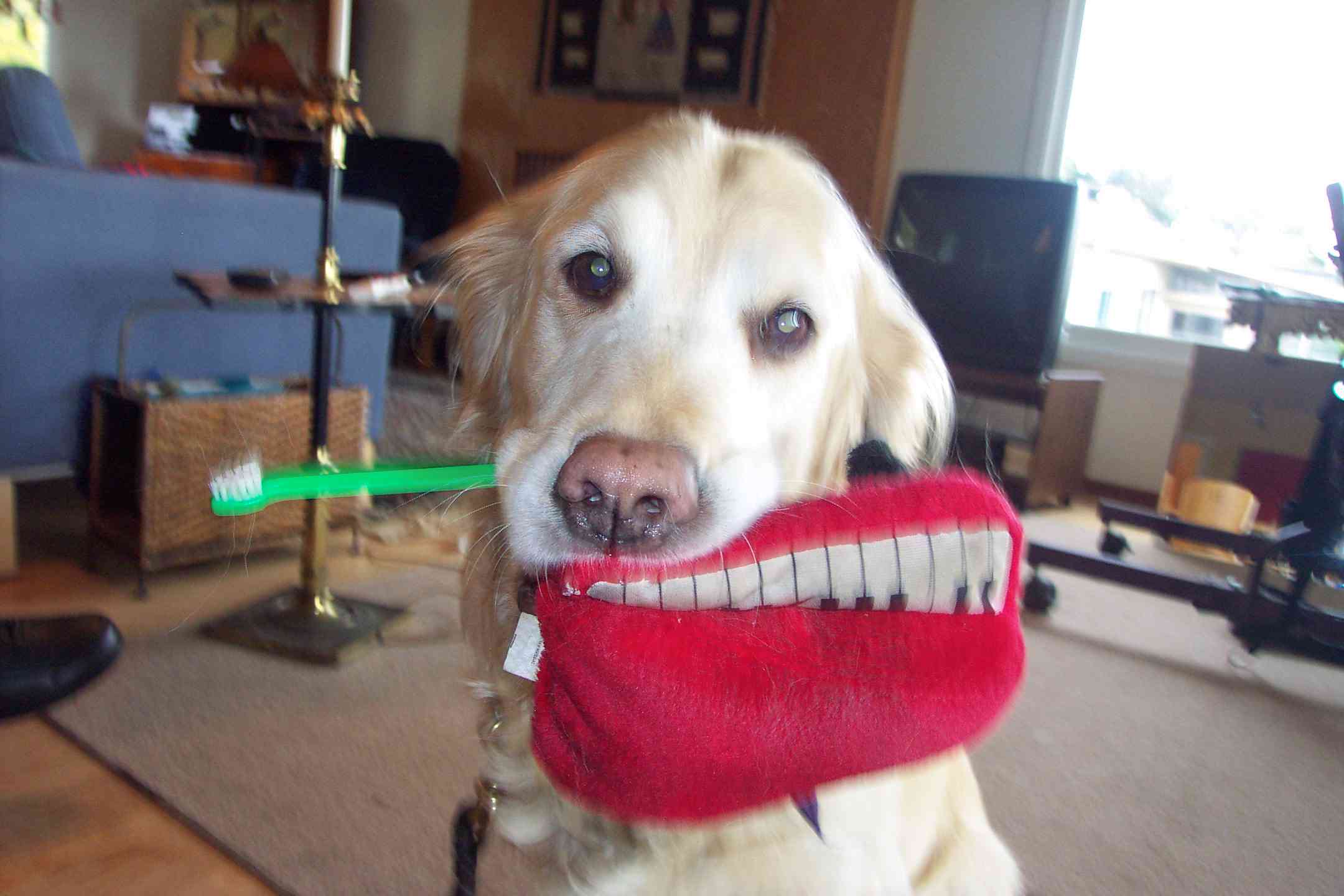}
\includegraphics[width=.4\textwidth]{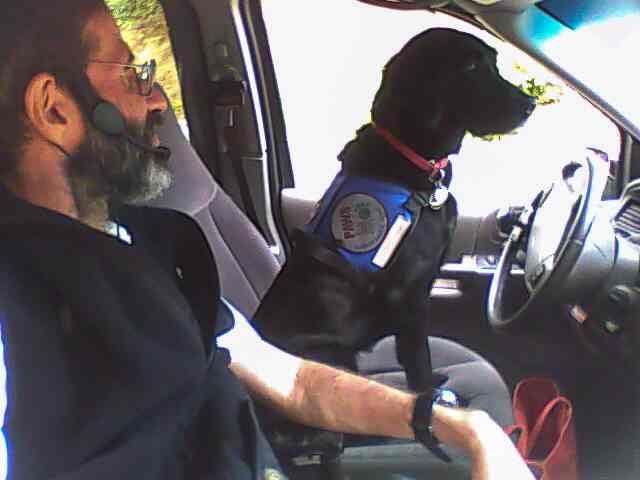}
\end{figure}

%   \hskip100pt \includegraphics[width=0.4\textwidth]{group.jpg}
  
% \begin{comment}
 Dan Sparks, a graduate student of Robert's, recalled Julep's second birthday party:
 \begin{quote} Several other dogs came to the party. The usual birthday stuff happened: Julep got presents and toys, as well as a nice big slice of cake. The cake was made special for the dogs from roast beef, with mashed potato icing.\end{quote}
 
 \noindent and went on to reflect:
 \begin{quote}
 % I am terribly sorry for his wife Tessa, and my condolences are with her. Still, it brings me unimaginable joy to know that the last years of his life were spent in that relationship, which I know made him very happy. This is something that I firmly believe Coleman deserved. You see, 
\ldots for all the amusing anecdotes (which reflect his unyielding sense of humor)\ldots 
% and mathematical ideas, 
perhaps the most memorable thing about him is how optimistic and positive he was\ldots
% And this despite his illness---it is truly, truly amazing. Whatever struggle I felt I had in life, or perhaps when I felt down or blue for whatever reason, it was never once that I still felt that way when leaving our weekly meetings. 
His positive energy was completely contagious, and the difficulties that he overcame put my own concerns into perspective. The net result of that was the greatest contribution that Robert Coleman made to my life.\end{quote}
% \end{comment}

%   \hskip100pt \includegraphics[width=0.4\textwidth]{Tessa}

%  His friend Marisa Castellano wrote:
% \begin{quote}Robert purchased an off-road wheelchair from [my husband] John back in 1989, and we used to go mountain biking with Robert and Bas [Edixhoven] back in the day when we still lived in El Cerrito. John and I would take turns towing Robert up the hills using bungee cords attached to our bikes. Sometimes we would tandem tow. And then Robert would fly downhill. Robert may have seemed meek in some ways, but he was a thrill seeker on those downhills!\end{quote}

Robert purchased an off-road wheelchair from John Castellano in 1989, which he used to go mountain biking with John, his wife Marisa, and Bas Edixhoven.  John and Marisa would take turns towing Robert up hills using bungee cords attached to their bikes.  According to Marisa, ``Sometimes we would tandem tow. And then Robert would fly downhill. Robert may have seemed meek in some ways, but he was a thrill seeker on those downhills!''  
 % Barry: I removed the following quote to shorten the piece a bit.  -- Matt
% As Kiran Kedlaya put it, ``[Robert's] ingenuity, humility, generosity, wry humor, and courage in the face of adversity (Banff in a wheelchair? In the middle of a snowy winter?!) are models for us all." 

\begin{figure}
\centering
\includegraphics[width=.4\textwidth]{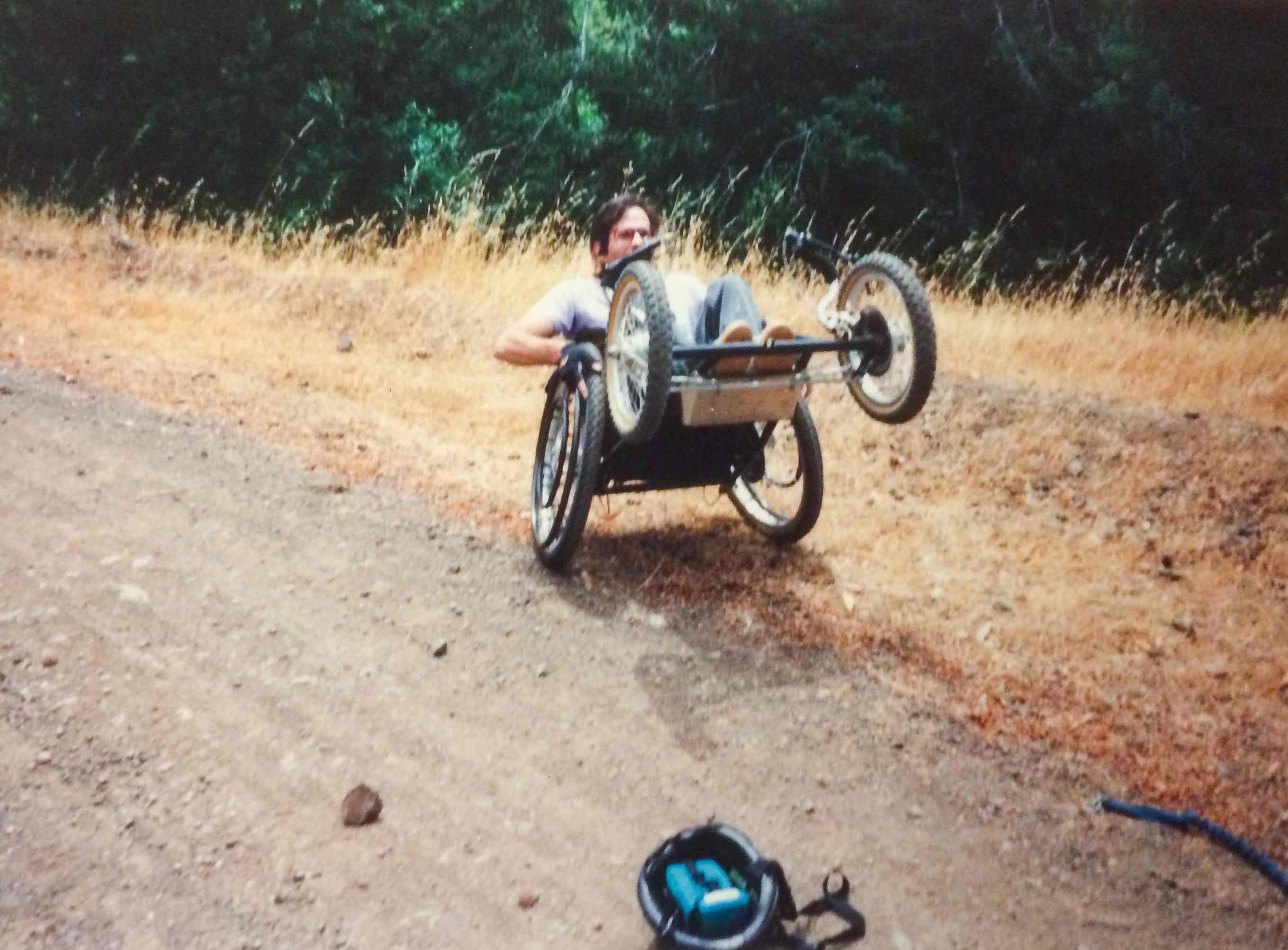}
\includegraphics[width=0.4\textwidth]{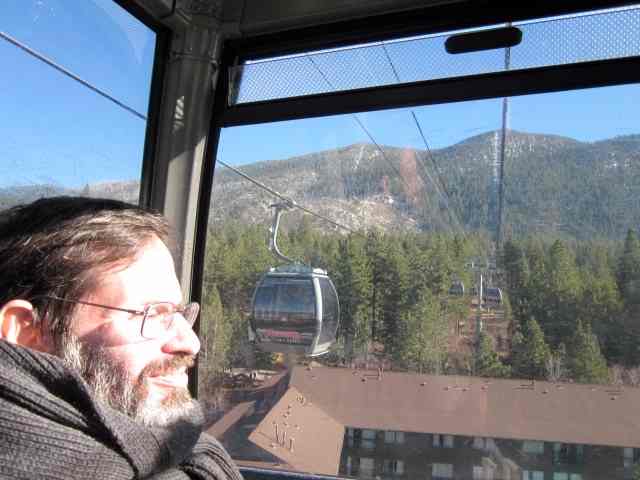}
\end{figure}
        
\section{Activism on behalf of people with disabilities}
        
% Robert worked tirelessly as a spokesperson for better planning and better services for the disabled\footnote{We thank Janet Abelson, Mayor of El % Cerrito, CA, for briefing us on this topic.}. He was a co-organizer of the Albany--El Cerrito access group,
% Committee for Disabled Rights in his own community of El Cerrito 
% which argued for---and succeeded in getting---a range of better resources for the disabled, including easier access to buildings and public 
% spaces.  For example, Robert succeeded in creating a long ramp for the disabled at an elementary school playground in El Cerrito: someone in a 
% wheelchair could ride the ramp to the top of a high slide and then slide down the ramp (without the wheelchair). 
% According to Janet Abelson, the mayor of El Cerrito and a co-organizer with Robert of the Albany--El Cerrito group, 
% Robert was a specialist in finding cases that could be turned into class action suits.  He worked closely with a local law firm that had filed several % such cases in the past, and succeeded multiple times in getting offending local businesses to cease their egregious behavior.  The money 
% obtained from such lawsuits was put to good use; for example, Robert's Albany--El Cerrito group helped install audible signals for the blind at 
% some major intersections in El Cerrito.  The group also paid for improvements to local libraries for patrons with impaired vision.
 
Almost immediately after the onset of his illness, Robert began
working as an activist and spokesperson for people with
disabilities\footnote{We thank Janet Abelson, Mayor of El Cerrito, CA,
for briefing us on this topic.}.  A resident of El Cerrito, CA, he
became an organizer of the Albany--El Cerrito Access Group, which
argued successfully for a range of better resources for the
disabled, including easier access to buildings and public spaces.
In particular, Robert managed to effect the installation of a long
ramp at an elementary school playground in El Cerrito that allowed
a wheelchair user to ride to the top of a high slide and then take
the slide back down to the ground (leaving the wheelchair behind).

\medskip

Robert seemed to have a special ability to identify situations
that could be turned into class action suits.  He worked closely
with a local law firm that had filed several such cases in the
past, and succeeded repeatedly in getting offending local
businesses to cease their egregious behavior.  The settlements from
these lawsuits were put to good use; for example, Robert's group
funded the installation of audible signals for the blind at major
intersections in El Cerrito.  The group also paid for improvements
to local libraries for patrons with impaired vision.

 \section{The Potato Society}   
 
For many years, Robert hosted a wine and cheese gathering in his office every Friday afternoon code-named ``Potatoes''.  Bill McCallum explains the origin of the name:
\begin{quote}
In Fall of 1978 or 1979 I purchased a case of Australian wine (Jacob's Creek Cabernet?) from the Wine and Cheese Cask and had it delivered to my cubicle in the math department at Harvard, meaning to take it home. Instead, a group of us started drinking it when the cheap Egri Bikaver the department served for Friday afternoon wine and cheese ran out\ldots 
% When the case was finished\ldots 
We started 
% going down to the Wine and Cheese Cask for supplies, 
having our own wine and cheese in my cubicle area. Greg Anderson, Robert Coleman, Robert Indik and Brad Osgood were regulars. In Fall of 1981 we moved into to Brad's office, and started going out to One Potato Two Potato for dinner afterwards. When Robert moved to Berkeley he continued the tradition there (it died out at Harvard after we all left).
\end{quote}
 
Robert's student Harvey Stein recalls:
 
\begin{quote}
Sometimes Potato Society meetings were a little more exotic. Hendrik Lenstra, returning from the Netherlands, once brought what must have been a gallon sized bucket of herring, along with a bottle of jenever. I spent much time filleting fish that night.
\end{quote}
 
And as Ken Ribet describes it:
 \begin{quote}
 We always had a great time at Potato Society. We could usually begin by opening up a couple of the bottles that were left over from the previous week. After Society members began to stream in, there were baguettes, new bottles of wine and chunks of cheese on the tables. We always managed to drag in extra people from the department who were walking near Robert's office. Our group often included postdocs, graduate students and staff members. There'd occasionally be an undergrad, but that was exceptional. We'd always have a good view of the Golden Gate Bridge and the setting sun. Depending on the season, the sun would be north of the bridge, right at the bridge or further south. Life stopped for a couple of hours late in the day every Friday. 
 % My wine glass from that era is still on one of my office bookshelves; it shows up sometimes in photos of office hours that are posted on my Facebook page.
 \end{quote}
 
% I removed the next quote to save space -- Matt
% There was a no-math rule at Potatoes, but this was broken on at least one occasion according to his (non-mathematician) friend John Castellano, who sold Robert his first off-road wheelchair:
%  \begin{quote}
% Robert was one of my first customers for my Off-Road Wheelchair design. We were neighbors, and started doing some rides together. He invited us for wine and cheese in his office, and that was our introduction to ``Potato Society''\ldots One time at Potatoes, I broke the anything-but-math rule and asked Robert about a design problem I was working on. I was designing a new style of bicycle frame with a shapely ``monocoque'' frame. I wanted to describe the perimeter mathematically so I could generate a 3-D model for stress analysis. I wanted to fit a curve going through six points on a plane. Then I could pick the endpoints, the slope near the endpoints, and a couple of points in between to control the shape\ldots
% [Robert] said ``it's just the summation over $i$ and $j = 1$ to $N$ of $y_i$ times $(x - x_j)/(x_i - x_j)$ for all $i \neq j$.'' I gulped. My (incorrect) expression took up a whole page, and his expression fit on half a line\ldots  I couldn't believe he was able to state it so compactly. Now after reading everybody's posts above, I see I was barely scratching the surface of Robert's understanding, but receiving the full measure of his generosity. Of course, it was not the only time.
% \end{quote}

\section{Coleman maps}  
   \label{reg}       
        
%  {\bf Barry and Ken: Insert something about Robert's Ph.D. thesis.}
 
Robert Coleman exploded onto the mathematical scene with his Ph.D. thesis, which was published as ``Division Values in Local Fields'' \cite{Coleman_DivisionValues} and continued in ``The Arithmetic of Lubin-Tate Division Towers'' \cite{Coleman_LubinTate}.  In the introduction to \cite{Coleman_DivisionValues}, Coleman writes:
\begin{quote}
In his work on cyclotomic fields, Kummer observed that various formal operations on power series had number theoretic applications. Perhaps the most striking of these was Kummer's idea of taking logarithmic derivatives of $p$-adic numbers. After a long period of neglect, various refinements and generalizations of Kummer's idea have recently been used by Iwasawa and Wiles to study explicit reciprocity laws, and by Coates and Wiles to study the arithmetic of elliptic curves with complex multiplication\ldots 
% Other applications of Kummer's observation include Iwasawa's explicit descriptions [I1], [[2] Of the Galois structure of various modules connected with local cyclotomic fields. 
The aim of the present paper is to begin a deeper and more systematic study of the local analytic theory which underlies these relations between power series and $p$-adic numbers.
\end{quote}

The main result of Coleman's paper is a general theorem on the interpolation of division values in Lubin--Tate formal groups which has found numerous applications to local class field theory, Iwasawa theory, and the arithmetic of $p$-adic $L$-functions.  We now describe this result, and its relation to Iwasawa theory from \cite{Coleman_LubinTate}, in the special setting of the (formal group associated to the) multiplicative group ${\mathbf G}_m$ following Pierre Colmez's exposition in \cite{Colmez_ColemanGeneralization}.  We also attempt to place these results into the larger context of $p$-adic representations and $p$-adic $L$-functions.

\medskip

To begin, % let $K$ be a finite extension of $\Q_p$ with valuation ring ${\mathcal O}_K$, and 
fix a compatible system $\{ \epsilon_n \}$ of roots of unity in $\overline{\Q_p}$ with $\epsilon_1 \neq 1$ and $\epsilon_{n+1}^p = \epsilon_n$ for all $n \geq 1$.  Let $K_n = \Q_p(\epsilon_n)$ and $K_\infty = \bigcup_{n \in {\mathbf N}} K_n$.  Let $\Gamma = {\rm Gal}(K_\infty / \Q_p)$ and let $\Lambda = {\mathbf Z}_p [[\Gamma]]$ be the completed group algebra of $\Gamma$. 
The {\em Coleman norm map} is a homomorphism $E: \varprojlim  {\mathcal O}_{K_n}^* \to \Lambda$ from the projective limit of the groups ${\mathcal O}_{K_n}^*$ with respect to the norm maps to the completed group algebra $\Lambda$.  To define the map $E$, Coleman first proves that given
an element $u = (u_n)$ of $\varprojlim  {\mathcal O}_{K_n}^*$, there is a unique invertible power series ${\rm Col}_u(T) \in \Z_p[[T]]$ such that
${\rm Col}_u(\epsilon_n - 1) = u_n$ for all $n \geq 1$.  (This is the interpolation theorem referred to above in the present context.)  The logarithmic derivative of 
${\rm Col}_u(T)$ again has coefficients in $\Z_p$ and there is a unique measure $\mu_u$ on $\Z_p$ such that
\[
\int_{\Z_p} (1 + T)^x \mu_u = (1+T) \frac{d}{dT} \log \left( {\rm Col}_u(T) \right).
\]
Restricting this measure to $\Z_p^*$ and pulling it back to $\Gamma$ via the cyclotomic character gives the Coleman map $E$, which turns out to be nearly an isomorphism.  As Colmez writes:

\begin{quote}
The measure which is used to define the Kubota--Leopoldt $p$-adic zeta function is the image of the cyclotomic units via this map, so the [Coleman norm map $E$] can be thought of as a machine producing $p$-adic $L$-functions out of compatible systems of units.
\end{quote}
 
 All of this can be interpreted in terms of the $p$-adic representation $\Q_p(1)$ associated to the cyclotomic character, and much later Cherbonnier and Colmez \cite{ColmezCherbonnier} used the theory of $(\phi,\Gamma)$-modules introduced by J.-M. Fontaine\footnote{Fontaine associates to any $p$-adic representation $V$ of the absolute Galois group of a local field $K$ of characteristic zero a so-called $(\phi,\Gamma)$-module, which is a certain module $D(V)$ equipped with a Frobenius operator $\phi$, a left inverse of $\phi$, and a continuous action of the $\Z_p$-cyclotomic Galois group of $K$ commuting with the action of $\phi$. This gives an equivalence between suitable categories of Galois representations and modules, and translates the study of $p$-adic representations V into the ``linear algebra'' of its corresponding module $D(V)$.} to extend Coleman's map to arbitrary $p$-adic representations.
 
 \medskip
 
Coleman subsequently realized that his norm map, combined with $p$-adic analytic properties of the dilogarithm function 
 \begin{equation}
 \label{eq:dilog}
 \ell_2(z) = \sum_{k=1}^\infty \frac{z^k}{k^2},
 \end{equation}
allow one to give a new explicit formula for the Hilbert norm residue symbol of local class field theory.
More specifically, in ``The Dilogarithm and the Norm Residue Symbol'' \cite{Coleman_NormResidueSymbol}, he presented a new formula for the norm residue symbol of exponent
$p^n$ attached to a cyclotomic extension of $\Q_p$ containing the $p^n$-th roots of unity.  This was the first concrete connection between the dilogarithm series, which plays an important role in Thurston's work on hyperbolic geometry and the work of Bloch and Milnor on algebraic K-theory, and class field theory.  (Such a connection is quite natural in view of the fact that in both settings the Steinberg relation between $x$ and $1-x$ plays a central role.)  Earlier explicit formulas for the norm residue symbol took on different forms for $p=2$ and odd primes, but in Coleman's formula all primes are treated equally.
Coleman relates being a local norm to being an exact differential, and is led by this to study a class of differential equations whose solutions can be expressed in terms of logarithms and dilogarithms.  He is also able to give a new proof of the non-degeneracy of the norm residue symbol from his explicit formula.

\medskip

Having been led by local class field theory to the study of p-adic analytic properties of the dilogarithm, one might naturally ask whether this function has a natural $p$-adic analytic continuation.
Over $\C$, the series given in (\ref{eq:dilog}) converges only in the open complex unit disc but can be analytically continued to a multi-valued function on $\C - \{ 0 \}$.  Since the power series in (\ref{eq:dilog}) also converges on the open $p$-adic unit disc in $\C_p$, Coleman had the idea to investigate whether the $p$-adic dilogarithm admits a similar natural extension.  The primary difficulty is that analytic continuation does not work, even in the context of Tate's rigid analytic spaces: the open disc is a maximal analytic domain for $\ell_2$.
% the fundamental theorems of rigid analysis introduced by John Tate, which were essentially the only results available for making sense of analytic continuation in the totally disconnected realm of the $p$-adics, do not seem very helpful in this context.  
Coleman's brilliant idea in ``Dilogarithms, Regulators, and $p$-adic $L$-functions'' \cite{Coleman_Dilogs} is to introduce what he called the ``Dwork principle'' of analytic continuation along Frobenius.\footnote{In the Acknowledgments to the paper, Coleman thanks Bernard Dwork ``who showed us that although one Frobenius may be best, all are good.''} Analytic extension for functions like polylogarithms, which satisfy differential equations with unipotent monodromy, is done by looking within the (vast) collection of locally analytic extensions at those global functions which satisfy certain rigidity assumptions on their twists by Frobenius.  Coleman shows that these Frobenius conditions allow one to prove uniqueness of the resulting extension.   

\medskip

Having extended the range of definition of the $p$-adic dilogarithm (and, more generally, all polylogarithms), Coleman derives a new formula---analogous to classical formulas for Dirichlet $L$-functions over the complex numbers---for the value of Kubota--Leopold $p$-adic $L$-functions at a positive integer $k \geq 2$ in terms of the $k^{\rm th}$ polylogarithm function.  He also uses the extended $p$-adic dilogarithm to produce a regulator map on $K_3(\C_p)$, analogous to Bloch's regulator map on $K_3(\C)$ in terms of the complex dilogarithm.  This uses a construction similar to Coleman's explicit formula for the Hilbert norm residue symbol, and a similar method gives a $p$-adic regulator map for Tate elliptic curves.  In the remarks at the end of the paper, Coleman writes: ``In a subsequent paper we intend to show how the ideas in this paper lead to a theory of $p$-adic abelian integrals.''  We turn to this theory next.

 \section{Coleman integration}
   \label{p-int}

The theory laid out in \cite{Coleman_Dilogs} can be viewed as a theory of $p$-adic Abelian integrals on ${\mathbf P}^1$.\footnote{The first such integral to appear in the literature was Kummer's $p$-adic logarithm, which he used in his work on explicit reciprocity laws and Fermat's Last Theorem.}  This approach, and its applications to $p$-adic $L$-functions, was extended to arbitrary curves with good reduction by Coleman and de Shalit in ``$p$-adic regulators on curves and special values of $p$-adic $L$-functions'' \cite{ColemandeShalit_Regulators}.  Specifically, if $C$ is a smooth complete curve over $\overline{\Q}_p$ whose Jacobian $J$ has good reduction, Coleman and de Shalit show how to define $p$-adic integrals for arbitrary meromorphic differential forms  on $C$.   As an application, if $W = H^0(C,\Omega^1_C)$ is the space of holomorphic differentials on $C$ and $T = {\rm Hom}(W,\overline{\Q}_p)$ is the tangent space of $J$ at the origin, they define a $p$-adic regulator pairing $r_{p,C} : K_2(\overline{\Q}_p(C)) \to T$ whose value at the Steinberg symbol $\{ f,g \}$ is the linear functional
\begin{equation}
\label{eq:padicreg}
r_{p,C}(\{ f,g \})(\omega) = \sum_{i=1}^t \int_{p_i}^{q_i} {\rm Log}(g) \cdot \omega,
\end{equation}
where ``Log'' denotes a fixed branch of the $p$-adic logarithm and ${\rm div}(f) = \sum_{i=1}^t (q_i) - (p_i).$

\medskip

In the special case where $C=E$ is an elliptic curve with complex multiplication, the pairing (\ref{eq:padicreg}) is related to a special value of the $p$-adic $L$-function of $E$; this gives a $p$-adic analogue of a theorem of Bloch.
As Coleman and de Shalit write in their paper:
\begin{quote}
The fragmentary evidence relating special values of classical $L$-functions to regulators on $K$-groups instigated very general and powerful conjectures of Bloch, Beilinson, and Deligne. Recently Soul{\'e} and Schneider have begun looking for $p$-adic conjectures. The relatively down-to-earth results of this paper, together with the earlier examples mentioned above, strongly support such analogues. Not less important, perhaps, is the indication that rigid analysis and $p$-adic integration ought to play some role in the proof of the $p$-adic conjectures.
\end{quote}

\begin{figure}[H]
\centering
\includegraphics[width=.4\textwidth]{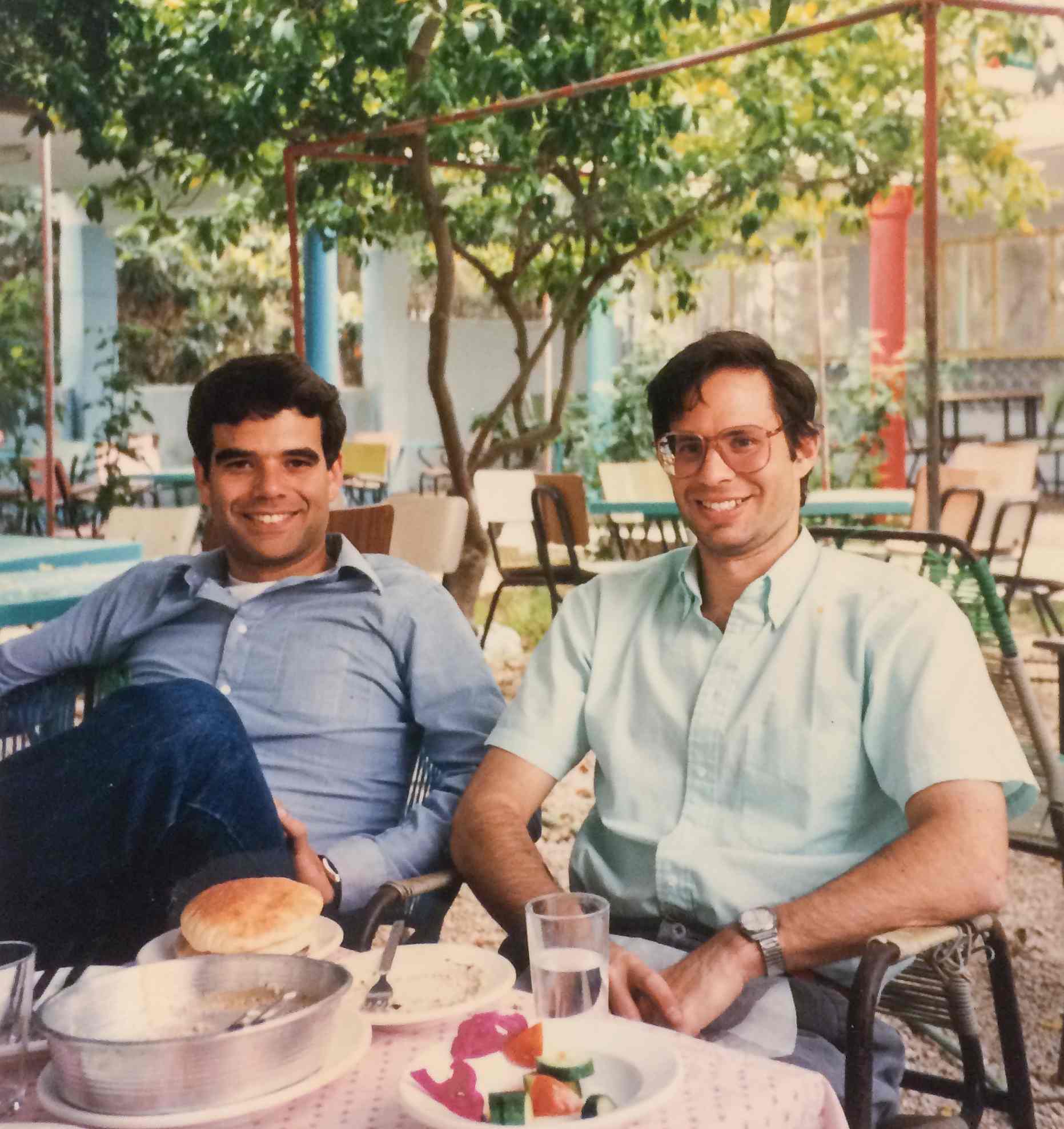}
% \caption{Robert and Udi de Shalit.}
\end{figure}

\medskip

For a more restrictive class of differential forms (differentials of the second kind\footnote{A meromorphic one-form on a variety $V$ is called a {\em differential of the second kind} if it is the sum of a holomorphic one-form and an exact one-form.}), but for varieties of any dimension having good reduction at $p$, Coleman developed in ``Torsion points on curves and $p$-adic Abelian integrals'' \cite{Coleman_AbelianIntegrals} a full-fledged theory of $p$-adic Abelian integrals of one-forms.  The construction is again based on the Dwork principle of ``analytic continuation along Frobenius''.  Coleman establishes the basic properties of these $p$-adic integrals, including an important functoriality result which he uses to show that $p$-adic Abelian integrals of holomorphic one-forms satisfy an addition law.  Once a branch of the logarithm has been fixed, the methods developed in \cite{Coleman_Dilogs,ColemandeShalit_Regulators,Coleman_AbelianIntegrals} can be combined to extend the theory of $p$-adic Abelian integrals to closed meromorphic one-forms on arbitrary varieties with good reduction, as well as to curves with bad reduction.  In the latter case, one does not always have a single-valued primitive for holomorphic differentials: there are periods which arise.  Coleman's theory of $p$-adic integration of one-forms was extended by Zarhin, Colmez, Besser, Vologodsky, and others (c.f. \cite{Breuil_SeminarBourbaki}), culminating in the recent book ``Integration of one-forms on $p$-adic analytic spaces'' by Vladimir Berkovich \cite{Berkovich_Integration}, which treats varieties of any dimension with no assumptions on the reduction type.  

\medskip

Going back to the paper \cite{Coleman_AbelianIntegrals}, as a consequence of his addition law Coleman shows that if a curve $C$ over $\C_p$ with good reduction at $p$ is embedded in its Jacobian $J$, the torsion points of $J(\C_p)$ lying on $C$ are the common zeros of certain $p$-adic Abelian integrals of the first kind on $C$.  These integrals can be expanded as power series in one variable (with respect to some local parameter on $C$), and the number of common zeros can be studied using Newton polygons and other classical tools from $p$-adic analysis.  Using this idea, in ``Ramified Torsion Points on Curves'' \cite{Coleman_RTPC} Coleman was able to give an entirely new proof of the Manin--Mumford conjecture, originally proved by Raynaud: if $C$ is a curve of genus at least 2 over a field of characteristic zero, then there are only finitely many torsion points of $J$ lying on $C$.  Unlike Raynaud's proof, Coleman's argument can sometimes be used to explicitly calculate this finite set of torsion points.  For example, together with Tamagawa and Tzermias, Coleman showed in ``The cuspidal torsion packet on the Fermat curve'' \cite{ColemanTamagawaTzermias} that if $C$ is the Fermat curve $X^m + Y^m + Z^m = 0$, embedded in its Jacobian $J$ via a {\em cusp} (i.e., a point $(x,y,z)$ with $xyz=0$), then the set of torsion points of $J$ lying on $C$ is precisely the set of cusps for all $m \geq 4$.  Coleman's theory of $p$-adic integration plays a crucial role in the proof.  Poonen \cite{Poonen_ComputingTorsionPoints}
presented an effective algorithm to calculate torsion points on curves by combining one of Coleman's main results in \cite{Coleman_RTPC}
with ideas of Buium from \cite{Buium_pjets}.

\medskip

Perhaps the biggest splash made by the theory of $p$-adic Abelian integrals was Coleman's highly influential paper ``Effective Chabauty'' \cite{Coleman_Chabauty}, in which he resurrected an old idea due to Chabauty and showed that it led to effective bounds for the number of rational points on an algebraic curve $C$ over a number field $K$, provided that the Mordell--Weil rank of the Jacobian of $C$ is not too large.  Specifically, Coleman proved the following theorem (which, for simplicity, we state only for $K={\mathbf Q}$):

\medskip

{\bf Theorem:} Let $C / {\mathbf Q}$ be an algebraic curve of genus $g \geq 2$, and assume that the rank $r$ of $J({\mathbf Q})$ is less than $g$.  Then for any prime $p > 2g$ at which $C$ has good reduction, 
\[
\# C({\mathbf Q}) \leq \# \overline{C}({\mathbf F}_p) + 2g-2.
\]

\medskip

Recall that the Mordell Conjecture, proved by Faltings (and independently by Vojta shortly thereafter), asserts that if $C$ is an algebraic curve over ${\mathbf Q}$ of genus at least 2, then the set $C({\mathbf Q})$ of rational points on $C$ is finite.  At the present time, however, we do not know an effective algorithm (even in theory) to compute this finite set. The techniques of Faltings and Vojta do lead in principle to an upper bound for the size of $C({\mathbf Q})$, but the bound obtained is very far from sharp and is hard to write down explicitly.   By contrast, when the Chabauty--Coleman method applies (i.e., when $r<g$), the bound is not only completely explicit, it is sometimes even sharp!\footnote{David Grant gives the example $y^2=x(x-1)(x-2)(x-5)(x-6)$ with $p=7$.} And even when Coleman's bound is not sharp, the method of proof can often be used (together with explicit computations) to find the set of rational points on $C$ exactly.  This has given rise to a whole industry in computational number theory, and many Diophantine equations have been solved using this method and its generalizations and refinements.  For example, in book 6 of the {\em Arithmetica}, Diophantus of Alexandria poses a question which comes down to finding the positive rational solutions to $y^2 = x^6 + x^2 + 1$, which describes a genus 2 curve $C$.
% (the only curve of genus greater than one which appears in the ten known books of the Arithmetica).  
Diophantus provides the solution $(1/2, 9/8)$ and essentially asks whether there are any other positive rational solutions.  In his 1998 Ph.D. thesis, Joe Wetherell used the method of Chabauty--Coleman to prove that there are no other solutions, thus resolving Diophantus' ancient problem.
% Nils Bruin used it in 1999 to determine all solutions in coprime integers to the equations $x^4+y^6=z^3$ and $x^2 + y^8 = z^3$
In 2007, Poonen, Schaefer, and Stoll \cite{PoonenSchaeferStoll} found all primitive integral solutions to the generalized Fermat equation $x^2 + y^3 = z^7$, combining the method of Chabauty--Coleman with the Frey--Serre--Ribet approach to Fermat's Last Theorem and the Wiles--Taylor modularity theorem.   The complete list of solutions with $xyz \neq 0$ is 
\[
% (\pm 1,-1,0),( \pm 1,0,1), \pm (0,1,1),
(\pm 3,-2,1),(\pm 71,-17,2),(\pm 2213459,1414,65),(\pm 15312283,9262,113),(\pm 21063928,-76271,17).
\]

\medskip

More recently, Michael Stoll developed a refinement of the method of Chabauty--Coleman, showing that there is a bound depending only on $g$ and $[K:Q]$ for the number of $K$-rational points on a hyperelliptic curve $C$ of genus $g$ over a number field $K$ with the property that the Mordell--Weil rank of its Jacobian is at most $g-3$.   And in another recent preprint, Poonen and Stoll prove that for $g \geq 3$, a positive fraction of hyperelliptic curves of odd degree $2g+1$ over ${\mathbf Q}$ have only one rational point, the point at infinity. They also prove a lower bound on this fraction that tends to 1 as the genus tends to infinity.  Their method combines a refinement of the Chabauty--Coleman method, based on an idea of McCallum, with the Bhargava--Gross equidistribution theorem for nonzero 2-Selmer group elements.  So the method of Chabauty--Coleman is alive and well.\footnote{Minhyong Kim has developed a partly conjectural non-abelian analogue of the method of Chabauty--Coleman which appears to hold great promise for the future; it is related to iterated Coleman integrals and the Grothendieck section conjecture.}

\medskip

To round off our discussion of Coleman integration, we mention that it has deep connections with $p$-adic Hodge theory which have presumably only begun to be understood.  One of the first such connections was made by Coleman himself in his 1984 paper
``Hodge--Tate periods and $p$-adic abelian integrals''  \cite{Coleman_HodgeTatePeriods}.  In that work, Coleman gives a new description of one of the maps involved in the Hodge--Tate decomposition of an abelian variety $A$ with good reduction.  Using this, he provides a formula for the Hodge--Tate periods\footnote{A curious feature of $p$-adic Hodge theory, which may be surprising to the uninitiated, is that the numbers which one routinely calls ``periods'' are almost never defined via integration of differential forms.} involving $p$-adic Abelian integrals on $A$.  
% A concrete consequence of this work is that a differential of the second kind on a curve with good reduction is exact if and only if its $p$-adic Abelian integral is bounded on the complement of a finite number of residue classes.  
Pierre Colmez built on this work in \cite{Colmez_Periods}, constructing $p$-adic periods for differentials of the second kind on arbitrary Abelian varieties $A$ over a local field $K$ with values in Fontaine's period ring $B^+_{\rm dR}$ via $p$-adic integration.

\medskip

Coleman continued to elucidate connections between $p$-adic integration and $p$-adic Hodge theory throughout his career.  For example, in the 1999 paper ``The Frobenius and monodromy operators for curves and Abelian varieties'' \cite{ColemanIovita_FrobeniusMonodromy}, written with Adrian Iovita, they show that if $A$ is an abelian variety over a local field $K$ having split semistable reduction, then 
% there is a compatibility between Fontaine's monodromy operator and Grothendieck's monodromy pairing.  
the $p$-adic integration pairing $T_p(A) \times H^1_{\rm dR}(A) \to B^+_{\rm dR}$
defined by Colmez in \cite{Colmez_Periods} induces a natural isomorphism between the integral part of $H^1_{\rm dR}(A)$ and Fontaine's ``monodromy module'' $D_{\rm st}(\Hom(V_p(A),\Q_p)$ which is compatible with the actions of the Frobenius and monodromy operators on both sides, as well as with the natural filtrations.  (Here $T_p(A)$ is the $p$-adic Tate module of $A$ and $V_p(A)=T_p(A) \otimes_{\Z_p} \Q_p$.)  As a concrete corollary, Coleman and Iovita obtain a proof of Fontaine's conjecture that if $A$ is an Abelian variety over a local field $A$, then $T_p(A)$ is crystalline if and only if $A$ has good reduction.

 \section{Coleman families}
   \label{eigen}

The classical work of Ramanujan on arithmetic properties of the Fourier coefficients of modular forms unearthed striking congruences which contain important number-theoretic information (see \cite{Ramanujan_Congruences}, as well as the discussion in \cite{Serre_Ramanujan}).  
% \cite{Serre_Ramanujan,Berndt-Ono_Ramanujan}).  
% See also \cite{Berndt-Ono_Ramanujan}
Perhaps the most famous  of these is his congruence modulo the prime $691$ relating
% the congruence modulo the prime $691$ that relates 
the Fourier coefficients of the cuspidal modular form  $$\Delta(q)\ =\ q\prod_{n=1}^{\infty}(1-q^n)^{24}\  =\ \sum_{n=1}^{\infty}\tau(n)q^n$$  to the Fourier coefficients of the Eisenstein series
 $$ E_{12}(q):= {\frac{1}{2}}\zeta(-11) \ + \ \sum_{n=1}^{\infty}\{\sum_{d \mid n}d^{11}\}q^n.$$ 
%  $$\tau(n) \ \equiv \   \sum_{d \mid n}d^{11} \ \pmod{691}$$  for every positive integer $n$ (where the sum is taken over positive divisors of $n$).
% (Here $q:=e^{2\pi i z}$).
Concretely, Ramanujan discovered that, for every positive integer $n$, $\tau(n)$ is congruent modulo 691 to the sum of the eleventh powers of all positive divisors of $n$.  
 
\medskip

The Fourier coefficients of the family of Eisenstein series of even integral weights $k =4,6,\dots,$ given explicitly by
 $$ E_{k}(q):= {\frac{1}{2}}\zeta(1-k) \ + \ \sum_{n=1}^{\infty}\{\sum_{d\mid n}d^{k-1}\}q^n,$$ 
already reveal a wealth of congruences between members of that family of different weights.  Indeed, for any odd prime $p$, if 
$k \equiv k'  \pmod{p-1}$ then Fermat's Little Theorem tells us that the non-constant Fourier coefficients of $E_{k}(q)$ are congruent mod $p$ to  the corresponding non-constant Fourier coefficients of $ E_{k'}(q)$, and the classical Kummer Congruence for Bernoulli numbers shows that the constant terms are congruent mod $p$ as well.  An appropriate elaboration of this observation allows us to construct a  family of 
% closely related 
``$p$-adic Eisenstein series" that vary over a $p$-adic analytic parameter space corresponding to their $p$-adic weights.
The Fourier coefficients in this family vary $p$-adically analytically as a function of the weight, and for any $p$-adic weight $k$ that happens to be a rational integer at least 2, the member of the family of that weight $k$ is an honest classical Eisenstein series of weight $k$.  Any two members of this family whose weights are appropriately close in this $p$-adic parameter space will have the property that their corresponding Fourier coefficients will be congruent modulo an appropriate power of $p$.  

\medskip

The work of Haruzo Hida (see e.g.~\cite{Hida_GaloisReps,Hida_OrdinaryHecke,Hida_Book}) was a major step forward in structurally organizing the ``wealth of congruences'' we've referred to.  Consider modular eigenforms on $\Gamma_0(pN)$ with $N$ not divisible by $p$. Such a modular form is said to be  of {\it slope zero}  if its eigenvalue for the Atkin--Lehner operator $U_p$ is a $p$-adic unit.  Hida showed that any such {\it cuspidal} eigenform fits into a $p$-adically varying family of  {\it $p$-adic} cuspidal modular eigenforms on $\Gamma_0(pN)$ (with nebentypus character depending on the weight) that have similar properties to the $p$-adic family of Eisenstein series described above. 
That is, the Fourier coefficients in Hida's families vary $p$-adically analytically over their parameter space, and for integral weights $k\ge 2$ the members of the family are classical cuspidal modular eigenforms. These {\it Hida families} are of finite degree over their parameter weight space, and provide congruences modulo powers of $p$ for $p$-adically nearby members.  Again, there are a wealth of such congruences.

\medskip
 
The existence of such ``$p$-adic continuity'' among cuspidal eigenforms might seem more surprising than in the Eisenstein family, 
in that Eisenstein series over the complex numbers fit into a {\it continuous spectrum} while cuspidal modular eigenforms are in the discrete spectrum in the classical harmonic analysis of the $L^2$ space of functions on the upper half plane relative to the action of an appropriate congruence subgroup of ${\rm SL}_2({\Z})$.  The $p$-adic setting, then, offers a type of coherence---via congruences modulo powers of $p$---that the classical one does not. 

\medskip

Robert Coleman's great contribution in \cite{Coleman_OverconvergentModularForms,Coleman_BanachSpaces} was to extend the construction of such $p$-adic families to cuspidal eigenforms on $\Gamma_0(pN)$
% on $\Gamma_0(p)$ 
that are merely assumed to be of finite slope; that is, not in the kernel of the $U_p$-operator.  
%  of arbitrary finite slope (meaning that the eigenvalue for the Atkin-Lehner operator $U_p$ is nonzero). 
{\it Coleman families} have all the properties---save one---that Hida families have.  That one difference, though, is crucial, and is why Coleman's contribution required a new vision, new techniques, and a strikingly original application of those techniques. Coleman families---in contrast with the Eisenstein families (which are of degree one over the parameter weight space) or Hida families (which are  algebraic, of finite degree over the weight space)---project only $p$-adically analytically to the parameter weight space.  Coleman obtains his families by performing a Fredholm analysis in a family of (infinite-dimensional) $p$-adic Banach spaces of $p$-adic modular forms (over $p$-adic weight space).  To achieve this was a real advance, and it led to the construction in \cite{ColemanMazur_Eigencurve} of the {\it eigencurve}, a curve which is (very likely) of infinite degree over the weight space that parametrizes all $p$-adic eigenforms of finite slope. The eigencurve is nowadays a basic tool for understanding the arithmetic of modular forms.
% a basic  ingredient in the inventory of tools used to  understand arithmetic problems related to modular forms. 
% There is now a great on-going project to construct {\it eigenvarieties} corresponding to automorphic representations for more general reductive groups.
   
\begin{comment}
%    \section{\bf \{A possible section on Robert's political work towards improving facilities for handicapped?\}}
%       \vskip20pt
%    \section{ \bf Coleman's theory of $p$-adic integration and his contributions to $p$-adic Hodge Theory}\label{p-int}  \section{\bf ``Coleman Families"}\label{eigen}  \section{\bf``Coleman maps"}\label{reg}
      \vskip20pt   \section{\bf \{A possible brief section on the 'potato society'?\}}
   \vskip20pt
      \section{} So many mathematicians have had their perspectives on life, and on thought, and on mathematics, deepened and transformed for the better because of their interaction and friendship with Robert.
\vskip20pt 
\centerline{{\bf Matt and Ken: }{\it We could elaborate with our personal reflections, or not, and continue with:}} 
 \vskip20pt   
 
\end{comment}

\section{Concluding remarks}
Robert Coleman was a kind, brave, and brilliant man whose influence on mathematics and on his friends and loved ones will long outlive his fragile body.  He focused on what is important in mathematics and in life, from beginning to end, and had the imagination, originality, and courage to work on problems which lay just at the edge of the possible.  He was unafraid to dream the great dreams of his subject.

\begin{figure}[h]
\centering
\includegraphics[width=.4\textwidth]{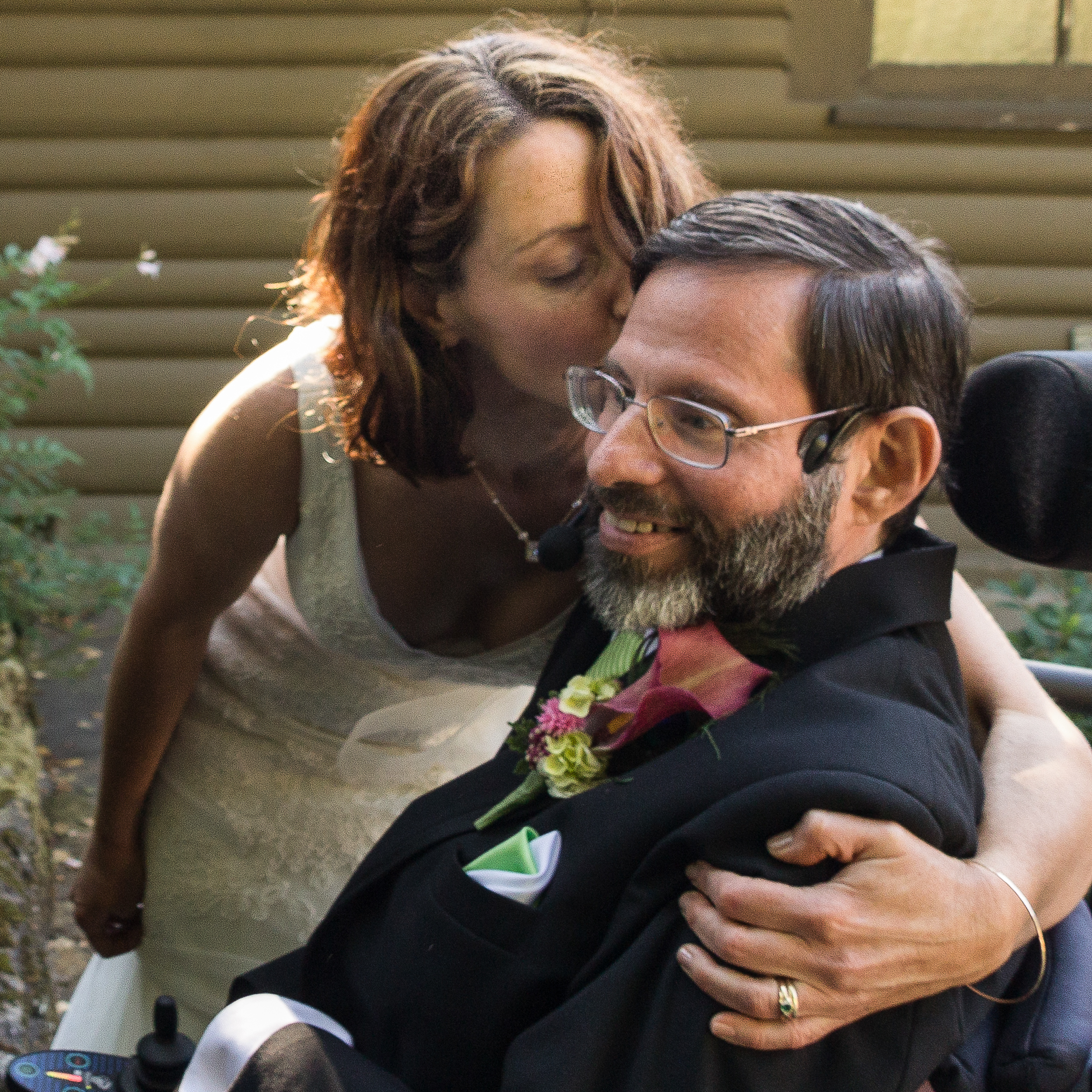}
\end{figure}

\medskip

A conference on $p$-adic Methods in Number Theory will take place in Berkeley, CA in May 2015 in honor of Robert's mathematical legacy. 

% \newpage

\bibliographystyle{alpha}
\bibliography{AMSNotices}

\end{document}